\documentclass[11pt]{article}

\usepackage[dvips,letterpaper,top=0.5in, bottom=0.5in, left=0.75in, right=0.5in,includefoot,heightrounded]{geometry}

\usepackage{srcltx}

\usepackage{amsmath}
\usepackage{amssymb,amsfonts}

\usepackage{abstract}

\usepackage{epsfig}
\usepackage[usenames,dvipsnames]{color}
\usepackage[usenames,dvipsnames]{xcolor}

\usepackage{booktabs}

\usepackage{setspace}
\usepackage{flushend}
\usepackage{url}

\usepackage[short,12hr]{datetime}

\usepackage{hyperref}

\usepackage[sc,tiny,center]{titlesec}
       \usepackage{mathptmx}

\newtheorem{proposition}{Proposition}

\newtheorem{lemma}{Lemma}

\newcommand{\fourier}{\stackrel{\mathcal{F}}{\longleftrightarrow}}
\newcommand{\cas}{\operatorname{cas}}
\newcommand{\hartley}{\stackrel{\mathcal{H}}{\longleftrightarrow}}
\renewcommand{\th}{${}^{\mathrm{th}}$}
\newcommand{\st}{${}^{\mathrm{st}}$}
\newcommand{\nd}{${}^{\mathrm{nd}}$}
\newcommand{\rd}{${}^{\mathrm{rd}}$}

\def\QED{\mbox{$\square$}}
\def\proof{\noindent{\it Proof:~}}
\def\endproof{\hspace*{\fill}~\QED\par\endtrivlist\unskip}



\title{%
A Factorization Scheme for Some Discrete Hartley Transform Matrices}

\author{%
H. M. de Oliveira%
\thanks{%
H. M. de Oliveira
was with 
Departamento de Eletr\^onica e Sistemas,
Universidade Federal de Pernambuco (UFPE).
He is currently
with
the Signal Processing Group,
Departamento de Estat\'{\i}stica, 
Universidade Federal de Pernambuco.
Email: \url{hmo@ufpe.br}
}
\and
R. J. Cintra%
\thanks{%
R. J. Cintra 
was with 
the Communications Research Group,
Departamento de Eletr\^onica e Sistemas,
Universidade Federal de Pernambuco.
He is currently
with
the Signal Processing Group,
Departamento de Estat\'{\i}stica, 
Universidade Federal de Pernambuco.
E-mail: \protect\url{rjdsc@de.ufpe.br}
}
\and
R. M. Campello de Souza%
\thanks{%
R. M. Campello de Souza
is with 
the Communications Research Group,
Departamento de Eletr\^onica e Sistemas,
Universidade Federal de Pernambuco.
Email: \url{ricardo@ufpe.br}
}
}

\date{\today\ @ \currenttime}
\date{}

\newcommand{\myabstract}{%
Discrete transforms such as the discrete Fourier transform 
(DFT) and the discrete Hartley transform (DHT) are
important tools in numerical analysis, signal processing, and
statistical methods.
The successful application of transform techniques relies on
the existence of efficient fast transforms.
In this paper some fast algorithms are derived.
The theoretical lower bound on the multiplicative complexity for the DFT/DHT
are achieved.
The approach is based on the factorization of DHT matrices.
Algorithms for short blocklengths such as 
$N \in \{3, 5, 6, 12, 24 \}$
are presented.
}

\newcommand{\mykeywords}{%
Discrete Hartley transform,
fast algorithms,
small blocklength
}

\begin{document}

  \maketitle
  \begin{abstract}
    \myabstract
  \end{abstract}
  \begin{center}
    \small
    \textbf{Keywords~\vspace{3mm}}
    \linebreak
    \mykeywords
  \end{center}
  \bigskip

\onehalfspacing

\section{Introduction}

Discrete transforms defined over finite or infinite fields
have been playing a relevant role in numerical analysis.
A striking example is the discrete Fourier Transform, which
has found applications in several areas.
Another relevant example concerns the discrete Hartley
transform (DHT)~\cite{Bracewell1}, the discrete version of the
integral transform introduced by Hartley in~\cite{Hartley}.
Besides its numerical side appropriateness,
the DHT has proven over the years to be a powerful 
tool~\cite{Bracewell2, Olejniczak, Shiu}. 
A decisive factor for applications of the DFT has
been the existence of fast transforms (FT) for
computing it~\cite{Blahut}. Fast Hartley transforms also exist and are
deeply connected to the DHT applications~\cite{Bi, Popovic}. Recent
promising applications of discrete transforms concern the
use of finite field Hartley transforms~\cite{Souza} to design digital
multiplex systems, efficient multiple access systems~\cite{Oliveira1}
and multilevel spread spectrum sequences~\cite{Oliveira2}.

Discrete transforms presenting a low multiplicative
complexity have been an object of interest for a long time. Very
efficient algorithms such as the Prime Factor Algorithm
(PFA) or Winograd Fourier Transform Algorithm (WFTA) have
also been used~\cite{Winograd, Yang}. The minimal multiplicative
complexity, $\mu(\cdot)$, of the one-dimensional DFT for all possible
sequence lengths, $N$, can be computed by converting the DFT
into a set of multi-dimensional cyclic convolutions. A lower
bound on the multiplicative complexity of a DFT is given 
in~\cite[Theorem 5.4, p.~98]{Heideman}. 
The values of $\mu_{DFT}(N)$ for some short
blocklengths are given in Table 1 (some local
minima of $\mu$).

The discrete Hartley transform of a signal $v_i$, $i=0,1,2,\ldots,N-1$ is defined by
\begin{equation}
V_k\triangleq
\sum_{i=0}^{N-1}
v_i \cas
\left(
\frac{2\pi ki}{N}
\right), \quad k=0,1,\ldots,N-1,
\end{equation}
where $\cas(x)=\cos(x)+\sin(x)$ is the ``cosine and sine''
Hartley symmetric kernel.

In this paper, some FTs are presented, which meet the
minimal multiplicative complexity. There is a simple
relationship between the DHT and the DFT of a given real
discrete signal $\mathbf{v}$. 
If $v_i \fourier  F_k$ is a DFT pair
and $v_i \hartley H_k$ is the corresponding DHT pair, 
then~\cite{Bracewell2} we have:
\begin{equation*}
H_k=\Re \big\{ F_k \big\} - \Im \big\{ F_k \big\}
\end{equation*}
and
\begin{equation*}
F_k=\frac{1}{2}
\Big[
(H_k+H_{N-k})-
j (H_k-H_{N-k})
\Big]
,
\end{equation*}
for $k=0,1,\ldots,N-1$.

\begin{table}%
\centering
\caption{Minimal multiplicative complexity for computing a DFT of length $N$}
\label{tabela}
\begin{tabular}{cc}
\toprule
$N$ & $\mu_{DFT}(N)$\\
\midrule
3&1\\
5&3\\
6&2\\
12&4\\
24&12\\
\bottomrule
\end{tabular}
\end{table}

Therefore, a fast algorithm for the DHT is also a fast algorithm for
the DFT and vice-versa~\cite[Corollary 6.9]{Heideman}. Besides
being a real transform, the DHT is also involutionary, i.e.,
the kernel of the inverse transform is the same as
the one of the direct transform (self-inverse transform).
Since the DHT is a more symmetrical version of a discrete
transform, this symmetry is exploited so as to derive a 
FT that requires the minimal number of real floating
point multiplications. 

\section{Computing the 3-point DHT}

Let
$\mathbf{v}=(v_0, v_1, v_2)^T\hartley\mathbf{V}=(V_0, V_1, V_2)^T$
be discrete Hartley transform pair of blocklength 3.
The matrix formulation of this transform corresponds to 
$\mathbf{V}=\mathbf{H}_3\mathbf{v}$,
where $\mathbf{H}_3$ is given by
\begin{equation}
\mathbf{H}_3=
\left [
\begin {array}{ccc}  
1& 1& 1\\
1& \frac{\sqrt 3 -1}{2}&- \frac{\sqrt 3 +1}{2}\\
1&- \frac{\sqrt 3 +1}{2}& \frac{\sqrt 3 -1}{2}
\end {array}
\right ].
\end{equation}

Note that the irrational elements of  $\mathbf{H}_3$
have the same decimal part, i.e., except from their
integer part, they have the same absolute value
(see that $(\sqrt 3 -1)/{2} \approx .366\!\ldots$ and
$(\sqrt 3 +1)/{2} \approx 1.366\!\ldots$).
So let us make the following decomposition:
\begin{equation}
\mathbf{H}_3=
\underbrace{
\left [
\begin {array}{ccc}  
1& 1& 1\\
1& \frac{\sqrt 3 -1}{2}&- \frac{\sqrt 3 -1}{2}\\
1&- \frac{\sqrt 3 -1}{2}& \frac{\sqrt 3 -1}{2}
\end {array}
\right ]
}_{\mathbf{H'}_3}
+
\left [
\begin {array}{ccc}
 & & \\
 & &-1\\
 &-1& 
\end {array}
\right ].
\end{equation}

Since the 2\nd\ and 3\rd\ columns of the new matrix $\mathbf{H'}_3$ have
the same elements (taken in absolute values),
we can consider new variables $v_1+v_2$ and $v_1-v_2$.
Thus, this substitution yields the
following matrix equation:
\begin{equation}
\left [
\begin {array}{c} 
v_0\\
v_1+v_2\\
v_1-v_2
\end {array}
\right ]
=
\left [\begin {array}{ccc} 
1&0&0\\
0&1&1\\
0&1&-1
\end {array}
\right ]
\left [
\begin {array}{c} 
v_0\\
v_1\\
v_2
\end {array}
\right ].
\end{equation}
So the transform can be expressed by:
\begin{equation}
\mathbf{V}=
\underbrace{
\left [
\begin{array}{ccc}
1& 1& \\
1& & \frac{\sqrt 3 -1}{2}\\
1& & -\frac{\sqrt 3 -1}{2}
\end{array}
\right ]
}_{\mathbf{H}''_3}
\left [
\begin {array}{c} 
v_0\\
v_1+v_2\\
v_1-v_2
\end {array}
\right ]
+
\left [
\begin{array}{ccc}
&&\\
&&-1\\
&-1&
\end{array}
\right]
\cdot
\mathbf{v}.
\end{equation}

Observe that the new matrix $\mathbf{H}''_3$ can be
splitted in two new matrices, as shown below:
\begin{equation}
\mathbf{H}''_3=
\left [
\begin {array}{ccc} 
1&1&\\
1&&1\\
1&&-1
\end {array}
\right ]
\cdot
\left [
\begin {array}{ccc} 
1&&\\
&1&\\
&&\frac{\sqrt 3 -1}{2}
\end {array}
\right ].
\end{equation}

Joining the above equations in a single
statement, we have that:
\begin{equation}
\begin{split}
\mathbf{V}=&
\left(
\underbrace{
\left [
\begin {array}{ccc} 
1&1&\\
1&&1\\
1&&-1
\end {array}
\right ]
}_{\mathbf{C}}
\cdot
\underbrace{
\left [
\begin {array}{ccc} 
1&&\\
&1&\\
&&a
\end {array}
\right ]
}_{\mathbf{B}}
\cdot
\underbrace{
\left [
\begin {array}{ccc} 
1&&\\
&1&1\\
&1&-1
\end {array}
\right ]
}_{\mathbf{A}}
+
\underbrace{
\left [
\begin {array}{ccc}  
&&\\
&&-1\\
&-1&
\end {array}
\right ]
}_{\mathbf{L}}
\right)
\mathbf{v}\\
\end{split},
\end{equation}
where $a=(\sqrt 3 -1)/{2}$.

One can recognize the pre-addition matrix $\mathbf{A}$,
the multiplication matrix $\mathbf{B}$ and
the post-addition matrix $\mathbf{C}$.
This algorithm introduces a new kind of matrices denoted by
$\mathbf{L}$. We will name them ``layer matrix''.

\begin{figure*}
\centering
\epsfig{file=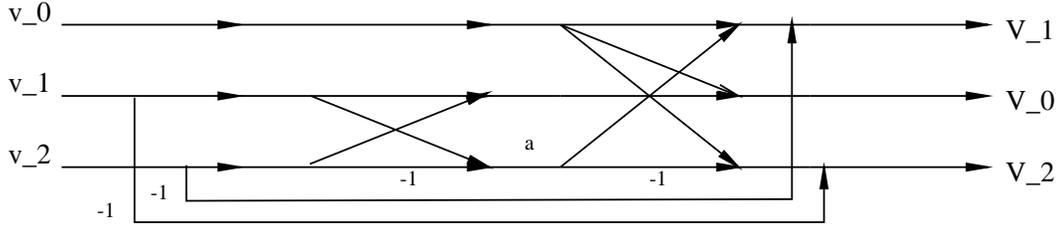, width=14cm}
\caption{3-point DHT fast algorithm diagram.}
\label{3-dht-o}
\end{figure*}

As the notation was explained, we can now express the entire
algorithm compactly by the following equation:
\begin{equation}
\mathbf{V=(CBA+L)v}
.
\end{equation}
This algorithms has only one nontrivial multiplication and 7 additions.
In Figure~\ref{3-dht-o}, there is a schematic diagram for
this transform.

\section{Computing the 5-point DHT}

Let $\mathbf{v} \hartley \mathbf{V}$ be a 5-point DHT pair.
The corresponding matrix formulation is now
$\mathbf{V}=\mathbf{H}_5 \mathbf{v}$, where
\begin{equation}
\mathbf{H}_5=
\left [
\begin {array}{ccccc} 
1&1&1&1&1\\
1&
a&%
b&%
c&%
d\\ %
1&
b&%
d&%
a&%
c\\ %
1&
c&%
a&%
d&%
b\\ %
1&
d&%
c&%
b&%
a%
\end {array}\right ],
\end{equation}
where
\begin{align*}
a&=1/4\left(\sqrt {5}-1+\sqrt {2}\sqrt {5+\sqrt {5}}\,\right),\\
b&=-1/4\left(\sqrt {5}+1-\sqrt {2}\sqrt {5-\sqrt {5}}\,\right),\\
c&=-1/4\left(\sqrt {5}+1+\sqrt {2}\sqrt {5-\sqrt {5}}\,\right),\\
d&=1/4\left(\sqrt {5}-1-\sqrt {2}\sqrt {5+\sqrt {5}}\,\right).
\end{align*}

We can combine the 2\nd\ and the 4\th\ columns as well as
the 3\rd\ and the 5\th\ ones using a Hadamard transform unit
of length~2 (a butterfly).
As a result of the process of combining columns, we achieve
the following matrix factorization:
\begin{equation}\label{matrix:5}
\mathbf{V}=
\mathbf{C}_3
\mathbf{C}_2
\mathbf{C}_1
\mathbf{B}
\mathbf{\tilde{A}}_2
\mathbf{A}_1
\mathbf{v},
\end{equation}
where $\mathbf{A}_1$ is the pre-addition matrix,
$\mathbf{B}$ is the multiplication matrix and
$\mathbf{C}_i$ are the post-addition matrices.
This matrices are detailed below:
\begin{alignat*}{2}
\mathbf{A}_1&=
\left[
\begin{smallmatrix}
1&&&&\\
&1&&&-1\\
&&1&-1&\\
&&1&-1&\\
&1&&&-1
\end{smallmatrix}
\right],
&\quad
\mathbf{\tilde{A}}_2&=
\left[
\begin{smallmatrix}
1&&&&\\
&1&1&&\\
&1&-1&&\\
&&&e&f\\
&&&f&-e
\end{smallmatrix}
\right],\\
\mathbf{B}&=
\left[
\begin{smallmatrix}
1&&&&\\
&\sqrt 5&&&\\
&&1&&\\
&&&1&\\
&&&&1
\end{smallmatrix}
\right],
&\quad
\mathbf{C}_1&=
\left[
\begin{smallmatrix}
1&1&&&\\
1&-1&&&\\
&&1&&\\
&&&1&\\
&&&&1
\end{smallmatrix}
\right],
\\
\mathbf{C}_2&=
\left[
\begin{smallmatrix}
1&&&&\\
&1&1&&\\
&1&-1&&\\
&&&1&\\
&&&&1
\end{smallmatrix}
\right],
&\quad
\mathbf{C}_3&=
\left[
\begin{smallmatrix}
1&&&&\\
&1&&1&\\
&&1&&1\\
&1&&-1&\\
&&1&&-1
\end{smallmatrix}
\right],
\end{alignat*}
where $e=\sqrt 2\sqrt{5-\sqrt 5}/2$ e $f=\sqrt 2\sqrt{5+\sqrt 5}/2$.

Now let us work in the matrix $\mathbf{\tilde{A}}_2$. 
Note that it contains multiplicative elements, namely $e$ and $f$,
and four additions.
We can go further and factorize this matrix in such a way that
purely multiplicative and additive matrices appear.
This can be done by the following method:
\begin{equation}
\label{decomp}
\mathbf{\tilde{A}}_2
=
\left[
\begin{array}{ccccc}
1&&&&\\
&1&1&&\\
&1&-&&\\
&&&1&\\
&&&&1
\end{array}
\right]
\left[
\begin{array}{cccccc}
1&&&&&\\
&1&&&&\\
&&1&&&\\
&&&1&1&\\
&&&&-1&1
\end{array}
\right]
\left[
\begin{array}{cccccc}
1&&&&&\\
&1&&&&\\
&&1&&&\\
&&&f+e&&\\
&&&&f&\\
&&&&&f-e
\end{array}
\right]
\left[
\begin{array}{ccccc}
1&&&&\\
&1&&&\\
&&1&&\\
&&&1&\\
&&&-1&1\\
&&&&1
\end{array}
\right].
\end{equation}
Thus, in Equation~\ref{matrix:5} one should
replace the matrix $\mathbf{\tilde{A}}_2$ by its
factorization.
The full decomposition of the original
transform matrix $\mathbf{H}_5$ is then achieved.

\begin{figure*}
\centering
\epsfig{file=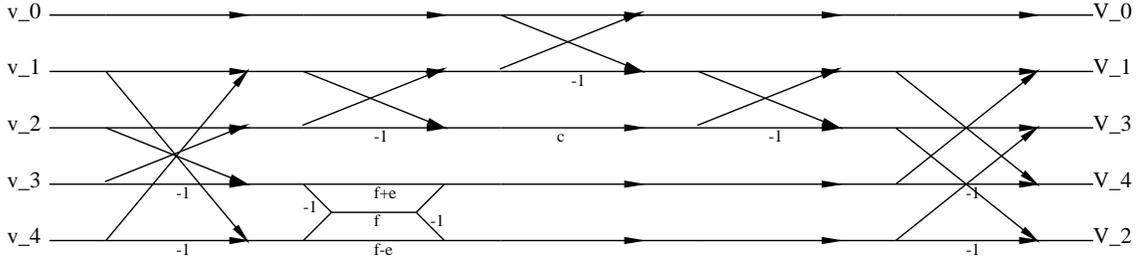, width=15cm}
\caption{5-point DHT fast algorithm.}
\label{5-dht-o-a}
\end{figure*}

The arithmetic complexity of this algorithm is 3
multiplications and 17 additions. The schematic diagram
is depicted in Figure~\ref{5-dht-o-a}.

\section{Computing a 6-point DHT}

Let us now consider $\mathbf{v} \hartley \mathbf{V}$ 
the transform pair related by the Hartley matrix
$\mathbf{H}_6$, where
\begin{equation}
\mathbf{H}_6=
\left [
\begin {array}{cccccc}
1&1&1&1&1&1\\
1& \frac{\sqrt 3 +1}{2}& \frac{\sqrt 3 -1}{2}&- 1&- \frac{\sqrt 3 +1}{2}&- \frac{\sqrt 3 -1}{2}\\
1& \frac{\sqrt 3 -1}{2}&- \frac{\sqrt 3 +1}{2}& 1& \frac{\sqrt 3 -1}{2}&- \frac{\sqrt 3 +1}{2}\\
1&- 1& 1&- 1& 1&- 1\\
1&- \frac{\sqrt 3 +1}{2}& \frac{\sqrt 3 -1}{2}& 1&- \frac{\sqrt 3 +1}{2}& \frac{\sqrt 3 -1}{2}\\
1&- \frac{\sqrt 3 -1}{2}&- \frac{\sqrt 3 +1}{2}&- 1& \frac{\sqrt 3 -1}{2}& \frac{\sqrt 3 +1}{2}
\end {array}
\right ].
\end{equation}

Using the Hadamard transform to combine
the 1\st\ and the 4\th\ columns,
the 2\nd\ and the 5\th\ columns,
and finally,
the 3\rd\ and the 6\th\ columns,
the matrix algorithm can reduced to:
\begin{equation}
\mathbf{V}=
\underbrace{
\left [
\begin {array}{cccccc} 
1& 1&1&&&\\
& & & 1& \frac{\sqrt 3 +1}{2}& \frac{\sqrt 3 -1}{2}\\
1& \frac{\sqrt 3 -1}{2}&- \frac{\sqrt 3 +1}{2}& & &\\
& & & 1& -1& 1\\
1&- \frac{\sqrt 3 +1}{2}& \frac{\sqrt 3 -1}{2}& && \\
& && 1& -\frac{\sqrt 3 -1}{2}& -\frac{\sqrt 3 +1}{2}
\end {array}\right ]
}_{\mathbf{H}'_6}
\mathbf{A}_1
\mathbf{v},
\end{equation}
where $\mathbf{A}_1$ is the first pre-addition matrix.
The matrix $\mathbf{A}_1$ is detailed below:
\begin{equation}
\mathbf{A}_1=
\left [\begin {array}{ccc|ccc} 
1&&&1&&\\
&1&&&1&\\
&&1&&&1\\\hline
1&&&-1&&\\
&1&&&-1&\\
&&1&&&-1
\end {array}\right ].
\end{equation}

The form of the matrix $\mathbf{A}_1$ is the same for
all transforms of even blocklength. The explanation to this
fact is given by the lemma below.
\begin{lemma}\label{A_1}
The pre-addition matrix
$\mathbf{A}_1$ 
of a Hadamard decomposition of a even blocklength DHT
has the following
construction:
\begin{equation}
\mathbf{A}_1=
\left[
\begin{array}{c|c}
\mathbf{I}_{\frac{N}{2}} & \mathbf{I}_{\frac{N}{2}}  \\ \hline
\mathbf{I}_{\frac{N}{2}} & -\mathbf{I}_{\frac{N}{2}}
\end{array}
\right]=
\mathbf{Had}_2 \otimes \mathbf{I}_{\frac{N}{2}},
\end{equation}
where $\mathbf{Had}_2$ is the Hadamard matrix,
$\otimes$ is the direct product and
$\mathbf{I}_{\frac{N}{2}}$ is a identity matrix
of order $N/2$.
\end{lemma}
\proof
The elements of the Hartley matrix, $\mathbf{H}_N$,
are governed by this property:
$h_{k,i+\frac{N}{2}} = (-1)^k h_{k,i}$, 
where $h_{k,i}$ is the
$(k,i)$-element of the transform matrix.
This property can be derived from
the $\cas(\cdot)$ arcs addition rule~\cite{Bracewell2}
$\cas(a-b)=\cos(b)\cas(a)-\sin(b)\cas'(a)$, 
where $\cas'(a)\triangleq\cos(a)-\sin(a)$.
Consequently we have that:
\begin{eqnarray*}
h_{k,i+\frac{N}{2}} &=&
\cas\left(\frac{2\pi k (i + \frac{N}{2})}{N}\right) \\
&=&\cas\left(\frac{2\pi ki}{N}+\pi k\right) \\
&=&(-1)^k\cas\left(\frac{2\pi ki}{N}\right) \\
&=& (-1)^k h_{k,i}
.
\end{eqnarray*}

Therefore the $i$\th\ and the $\left(i+\frac{N}{2}\right)$\th\ columns
have the same absolute value, which allows us to
combine them via Hadamard transform.
New variables arise from this technique:
$(v_i+v_{i+\frac{N}{2}})$ and $(v_i-v_{i+\frac{N}{2}})$, 
$i=0,\ldots, N/2-1$.
These new variables are generated by the matrix $\mathbf{A}_1$.
\endproof

Using the same strategy described in the 3-point DHT algorithm,
we can take aside the integer part of some elements of the
matrix 
$\mathbf{H}'_6$.
This procedure yields to a new more ``balanced'' matrix.
These steps are represented by the following equation:

\begin{equation}
\begin{split}
\mathbf{H}'_6=&
\left [
\begin {array}{cccccc} 
1& 1&1&&&\\
& & & 1& \frac{\sqrt 3 +1}{2}& \frac{\sqrt 3 -1}{2}\\
1& \frac{\sqrt 3 -1}{2}&- \frac{\sqrt 3 +1}{2}& & &\\
& & & 1& -1& 1\\
1&- \frac{\sqrt 3 +1}{2}& \frac{\sqrt 3 -1}{2}& && \\
& && 1& -\frac{\sqrt 3 -1}{2}& -\frac{\sqrt 3 +1}{2}
\end {array}\right ] \nonumber \\
=&
\underbrace{
\left [
\begin {array}{cccccc} 
1& 1&1&&&\\
& & & 1& \frac{\sqrt 3 -1}{2}& \frac{\sqrt 3 -1}{2}\\
1& \frac{\sqrt 3 -1}{2}&- \frac{\sqrt 3 -1}{2}& & &\\
& & & 1& -1& 1\\
1&- \frac{\sqrt 3 -1}{2}& \frac{\sqrt 3 -1}{2}& && \\
& && 1& -\frac{\sqrt 3 -1}{2}& -\frac{\sqrt 3 -1}{2}
\end {array}\right ]
}_{\mathbf{H}''_6}+
\underbrace{
\left [
\begin {array}{cccccc} 
&&&&&\\
&&&&1&\\
&&-1&&&\\
&&&&&\\
&-1&&&&\\
&&&&&-1
\end {array}
\right ]
}_{\mathbf{L}}.\\
\end{split}
\end{equation}
See that
$\mathbf{V}=
\left(
\mathbf{H}''_6
+
\mathbf{L}
\right)
\mathbf{A}_1
\mathbf{v}$.

Carrying out the procedure of combining columns which ``agree'',
we will have the next pre-addition matrix $\mathbf{A}_2$:
\begin{equation}
\mathbf{A}_2=
\left [
\begin {array}{cccccc} 
1&&&&&\\
&1&1&&&\\
&1&-1&&&\\
&&&1&&\\
&&&&1&1\\
&&&&1&-1
\end {array}
\right ].
\end{equation}
This makes the matrix  $\mathbf{H}'_6$ be written as
$\mathbf{H}'_6=\mathbf{H}''_6\cdot\mathbf{A}_2$,
as seen in this equation:
\begin{equation}
\mathbf{H}''_6=
\underbrace{
\left [
\begin {array}{cccccc} 
1& 1&&&&\\
& & & 1& \frac{\sqrt 3 -1}{2}& \\
1& & \frac{\sqrt 3 -1}{2}& & &\\
& & & 1& & -1\\
1&& -\frac{\sqrt 3 -1}{2}& && \\
& && 1& -\frac{\sqrt 3 -1}{2}& 
\end {array}\right ]
}_{\mathbf{H}'''_6}
\mathbf{A}_2.
\end{equation}

Now observe that
the factorization of 
$\mathbf{H}'''_6$
yields the multiplication matrix,
$\mathbf{B}$, and the post-additions
matrix,
$\mathbf{C}$.
\begin{equation}
\mathbf{H}'''_6=
\underbrace{
\left [
\begin {array}{cccccc} 
1& 1&&&&\\
& & & 1& 1& \\
1& & 1& & &\\
& & & 1& & -1\\
1&& -1& && \\
& && 1& -1& 
\end {array}\right ]
}_{\mathbf{C}}
\underbrace{
\left [
\begin {array}{cccccc} 
1&&&&&\\
&1&&&&\\
&&a&&&\\
&&&1&&\\
&&&&a&\\
&&&&&1
\end {array}\right ]
}_{\mathbf{B}}
,
\end{equation}
where $a=(\sqrt 3 -1)/{2}$.

\begin{figure*}
\centering
\epsfig{file=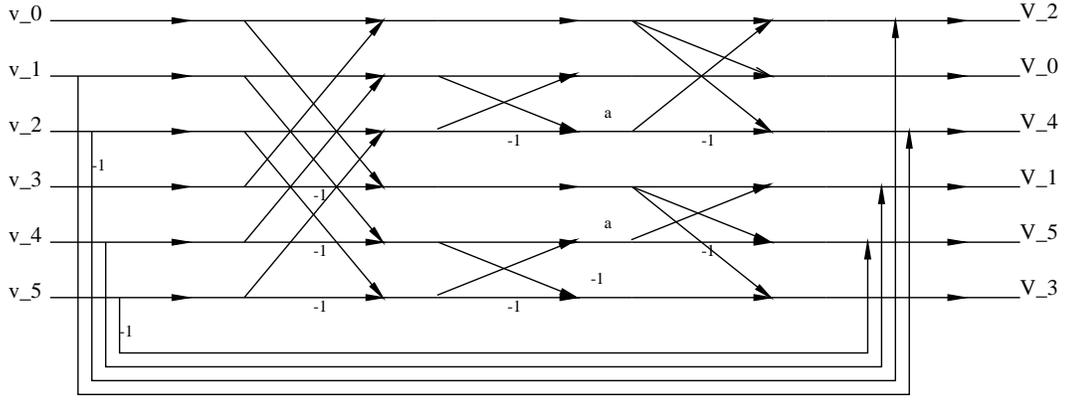, width=14cm}
\caption{6-point DHT fast algorithm diagram.}
\label{6-dht-o}
\end{figure*}

We have then completed the algorithm, and it can be represented by:
\begin{equation}
\mathbf{V}=(\mathbf{C}\mathbf{B}\mathbf{A}_2+\mathbf{L})\mathbf{A}_1\mathbf{v}.
\end{equation}
This algorithm has two multiplications and 20~additions and
is depicted in Figure~\ref{6-dht-o}.

\section{Computing the 12- and the 24-point DHT}

The procedure used to derive the 3- and 6-point DHT fast transforms
can be extended to other blocklengths, such as 12 and~24.
We derived these algorithms and achieved 
the  arithmetic complexity showed in Table~\ref{12:24}.

\begin{table}
\centering
\caption{Arithmetic complexity for the proposed 12- and 24-point DHT fast algorithm.
The function $\alpha(N)$ returns the additive complexity of the
implementation.}
\label{12:24}
\begin{tabular}{ccc}
\toprule
$N$ & $\mu(N)$ & $\alpha(N)$\\
\midrule
12&4&52\\
24&12&138\\
\bottomrule
\end{tabular}

\end{table}

\begin{figure*}
\centering
\epsfig{file=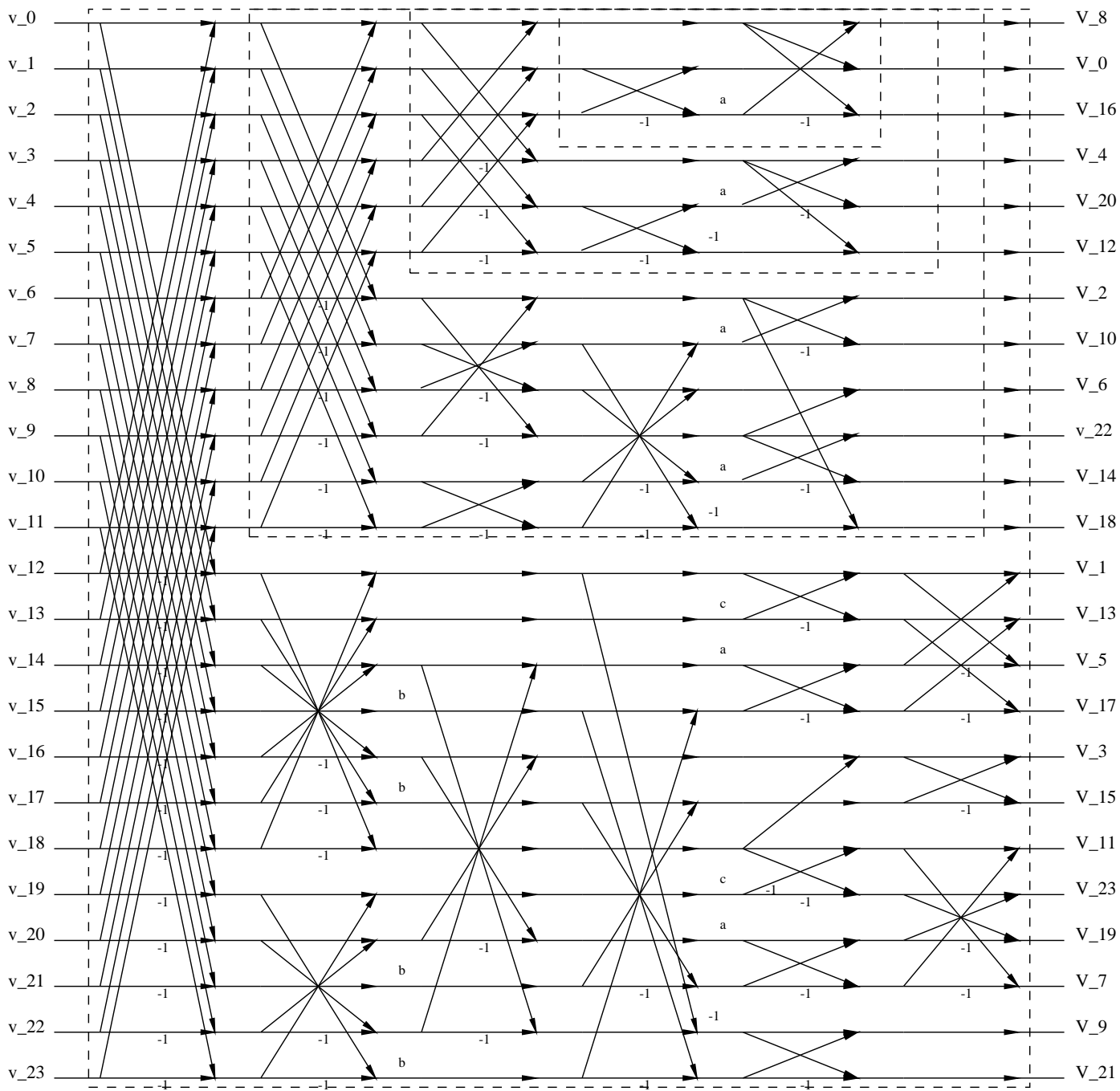, width=\linewidth}
\caption{24-point DHT fast algorithm diagram (derivations omitted). Shorter transforms
are embedded.}
\label{all-in-1}
\end{figure*}

In Figure~\ref{all-in-1}, we see a diagram of the 24-point DHT fast
transform, where the shorter transforms (3-, 6-, 12-point)
are embedded.

The algorithms proposed so far can be described in a
general framework according to the following proposition.
\begin{proposition}
The DHT decomposition has the following general formulation
\begin{equation}
\begin{split}
\mathbf{V}=&
\bigg(
  \Big(
    \big(
      (
        \mathbf{C}_n\mathbf{B}_n\mathbf{A}_n+\mathbf{L}_{n-1}
      )
      \mathbf{C}_{n-1}\mathbf{B}_{n-1}\mathbf{A}_{n-1} \cdots + \mathbf{L}_2
    \big)
    \mathbf{C}_2\mathbf{B}_2\mathbf{A}_2+\mathbf{L}_1
  \Big)
  \mathbf{C}_1\mathbf{B}_1\mathbf{A}_1+\mathbf{L}_0
\bigg)
\mathbf{v},\\
\end{split}
\end{equation}
where $n$ is the number of ``layers'' in the decomposition.\endproof
\end{proposition}

\section{Conclusions}

Short blocklength DHT fast algorithms that achieve the lower
bound on the multiplicative complexity were derived.
Low values for additive complexity were also found.
These algorithms can be implemented
in
digital signal processors capable
of low-power consumption.

\section*{Acknowledgments}

This work was partially supported by CNPq and CAPES.

{\footnotesize
\bibliographystyle{IEEEtran}
\bibliography{ref}
}

\end{document}